\documentclass[10pt]{amsart}

\usepackage{amsmath}
\usepackage{hyperref}
\usepackage{graphicx}
\usepackage{epsfig}
\usepackage{latexsym}
\usepackage{bm}
\usepackage{amssymb}

\def\diag{\mathop{\rm diag}}

\def\TV{{\rm{TV}}}

\def\diag{\mathop{\rm diag}}

\def\diag{\mathop{\rm diag}}
\def\sign{\mathop{\rm sign}}

\theoremstyle{remark}

\numberwithin{equation}{section}

\usepackage{hyperref}

\begin{document}

\title[On a nonlinear Laplacian based filter for noise removal]{On a nonlinear Laplacian based filter for noise removal}

\author{N. S. Hoang}

\address{Department of Computing and Mathematics, University of West Georgia, Carrollton, GA 30118, USA}

\email{nhoang@westga.edu}


\date{}


\begin{abstract}
We propose a nonlinear filter for noise removal based on the Laplacian for 1D and 2D data. The method utilizes the solution to a fourth-order nonlinear PDE involving the Laplacian for data reconstruction. Evolution equations are introduced to solve this fourth-order nonlinear equation. Numerical experiments show that the new filter preserves discontinuities while filtering out noise. The restored data are piecewise linear and avoid the staircase effect commonly observed with total variation denoising methods.
\end{abstract}
\maketitle

\section{Introduction}
Removing noise from noisy data is crucial in various applications and has been extensively explored in the literature. Linear filters, such as the Gaussian filter, often fail to preserve edges and discontinuities effectively. One of the most popular denoising methods is Total Variation (TV) denoising, introduced by Rudin et al. \cite{ROF}. Given a function $u_0$ defined on a bounded domain 
$D$ in $\mathbb{R}^2$, the TV denoising method seeks to minimize the following problem:
\begin{equation}
\label{eq1}
\min_{u}\int_D |\nabla u| dx
\end{equation}
with constraints
$$
\int_D u\,dx = \int_D u_0\,dx,\qquad \int_D (u - u_0)^2\,dx = \eta^2
$$
as the restored data of the noisy function $u_0$. 
Instead of directly solving the complicated problem \eqref{eq1}, Rudin et al. \cite{ROF} suggested solving the following unconstraint minimization problem:
\begin{equation}
\label{eq2}
\min_{u} \TV(u) + \frac{\lambda}{2} \|u - u_0\|^2,\qquad \TV(u) := \int_D |\nabla u|\, dx. 
\end{equation}
One of the approaches to solve problem \eqref{eq2} is to solve the 
Euler-Lagrange equation associated to \eqref{eq2} as follows
\begin{align}
\label{x1}
0 &= \nabla \cdot \frac{\nabla u}{|\nabla u|} - \lambda (u - u_0)\quad \text{in}\quad D,\\
\label{x2}
\frac{\partial u}{\partial n} &= 0\quad \text{on}\quad \partial D.
\end{align}
To solve equations \eqref{x1} and \eqref{x2}, it was suggested in \cite{ROF} to solve the following 
evolution equation
\begin{align}
\label{x3}
u_t &= \nabla \cdot \frac{\nabla u}{|\nabla u|} - \lambda (u - u_0)\quad \text{in}\quad D,\\
\label{x4}
u(x,0)&= \hat{u}_0,\qquad \frac{\partial u}{\partial n} = 0\quad \text{on}\quad \partial D,
\end{align}
and use the solution $u(t,x)$ for large $t$ as the reconstruction of $u_0$. 
For appropriate values of $\lambda$, the solution $u(t,x)$ to equations \eqref{x3} and \eqref{x4} converges to the solution ofl equations \eqref{x1} and \eqref{x2}.

Since its introduction, the TV algorithm has become a significant tool in signal processing. However, its main drawback is that the reconstructed functions are often piecewise constant, resembling a staircase: a phenomenon commonly referred to as the `staircase effect' in TV denoising. 
Numerous studies have focused on improving the TV model for denoising ( see, e.g., \cite{BKP}, \cite{BCM},\cite{CMM}, \cite{ZC}, \cite{SC}, \cite{STC}, \cite{GQSW}). Among these efforts, the total variation $TV(u)$ in \eqref{eq2} has been replaced by $\int_D\varphi(|\nabla u|)\, dx$ or  $\int_Dg(x)|\nabla u|\, dx$ (see, e.g., \cite{CMM}).

In \cite{CMM} Chan et al. the following model was proposed and studied
\begin{align}
\label{eqj1}
\int_D \big[\alpha |\nabla u|_\beta + \mu \Phi(|\nabla u|) (\mathcal{L}(u))^2\big]\, dx  +  \frac{1}{2} \|u - u_0\|^2_{\mathcal{L}^2} \longrightarrow \min_u,
\end{align}
where
$$
|x|_\beta = \sqrt{x^2 + \beta},\qquad \Phi(x) = 1/(\sqrt{x^2 + \gamma})^p,\qquad p>0.
$$

Numerical experiments in \cite{CMM} have demonstrated improvements over the standard Total Variation method by Rudin et al. (see  \cite{CMM}). 
However, the equation used to solve problem \eqref{eqj1} is quite complex. 
This paper aims to introduce a simpler model than the one in \eqref{eqj1}, which effectively avoids the `staircase effect' in TV denoising.

\section{A nonlinear filter}

\subsection{1-D noise removal}

Let $u_0$ be a noisy real-valued function. We proposed using $u$, the solution to the equation below, as the restored data for $u_0$:
\begin{gather}
\label{25a1}
0 = \frac{d^2}{dx^2}\bigg(\frac{u''}{\big[(u'')^2 + \epsilon\big]^p} \bigg) + \lambda(u - u_0)\quad in \quad D = [a,b],\qquad u'':=\frac{d^2u}{dx^2},\\
\label{25a2}
\frac{du}{dx}\bigg|_{x = a} = \frac{du}{dx}\bigg|_{x = b} = \frac{d^2u}{dx^2}\bigg|_{x = a} = \frac{d^2u}{dx^2}\bigg|_{x = b} =0.
\end{gather}
Here, $\lambda >0$ and $p\ge 0.5$ are parameters to be tuned experimentally. The constant $\epsilon>0$ is introduced to prevent division by zero errors. The parameter $\lambda$ typically depends on the noise level. 

Equation \eqref{25a1} with boundary condition \eqref{25a2} is a nonlinear forth-order equation. Since our primary goal is to test the new model (rather than to solve equation \eqref{25a2} efficiently), we choose to address equation \eqref{25a1} with boundary condition \eqref{25a2} by solving the following evolution equation:
\begin{gather}
\label{25a3}
\frac{du}{dt}= -\frac{d^2}{dx^2}\bigg(\frac{u_{xx}}{\big[(u_{xx})^2 + \epsilon\big]^p} \bigg) - \lambda(u - u_0)\qquad in \qquad D = [a,b],\\
\label{25a4}
\frac{du}{dx}\bigg|_{x = a} = \frac{du}{dx}\bigg|_{x = b} = \frac{d^2u}{dx^2}\bigg|_{x = a} = \frac{d^2u}{dx^2}\bigg|_{x = b} =0,\qquad u(0,x) = u_0(x),\quad x\in D,
\end{gather}
and use $u(t)$, for sufficiently large $t>0$, as the restored data for $u_0$. 

If the noise level $\|u_0 - f\| = \delta$ is known, the parameter $\lambda$ in equation \eqref{25a3} 
can be selected using a method similar to the one in \cite{ROF} as follows. At the equilibrium solution we have $\frac{du}{dt} = 0$. From equation \eqref{25a3} one obtains
$$
-\frac{d^2}{dx^2}\bigg(\frac{u_{xx}}{\big[(u_{xx})^2 + \epsilon\big]^p} \bigg) = \lambda(u - u_0). 
$$
This implies
$$
-(u - u_0)\frac{d^2}{dx^2}\bigg(\frac{u_{xx}}{\big[(u_{xx})^2 + \epsilon\big]^p} \bigg) = \lambda(u - u_0)^2.
$$
Integrate this equation over the interval $[a,b]$ we get
\begin{equation}
\label{eqlambda}
- \int_a^b (u - u_0)\frac{d^2}{dx^2}\bigg(\frac{u_{xx}}{\big[(u_{xx})^2 + \epsilon\big]^p} \bigg)\, dx = \lambda\int_a^b (u - u_0)^2\, dx = \lambda \delta.
\end{equation}
The parameter $\lambda$ is then determined from equation \eqref{eqlambda}.

Assume that $u(x)$ is available at equidistance nodes $x_i = x_0 +ih$, $i = 1,...,N$, $x_0:=a$, $h = (b-a)/N$. We use the following formula to compute the second derivatives $u''(x_i)$ 
\begin{equation}
\label{4m1}
u''(x_i) = \frac{u(x_{i-1}) - 2u(x_i) + u(x_{i+1})}{h^2} + O(h^2),\qquad i=1,...,N-1.
\end{equation}
The boundary condition $u'(a) = u'(b)$ becomes:
\begin{equation}
\label{4m2}
u(x_0) = u(x_1),\qquad u(x_{N-1}) = u(x_N).
\end{equation}
This implies
\begin{equation}
\label{4m4}
u''(x_1) \approx \frac{u(x_{0}) - 2u(x_1) + u(x_{2})}{h^2} =  \frac{- u(x_1) + u(x_{2})}{h^2}
\end{equation}
and 
\begin{equation}
\label{4m5}
u''(x_{N-1}) \approx \frac{u(x_{N-2}) - 2u(x_{N-1}) + u(x_{N})}{h^2} =  \frac{- u(x_{N-1}) + u(x_{N})}{h^2}.
\end{equation}

The condition $u''(a) = u''(b)=0$ can be satisfied by defining the values of $u$ at $x_{-1} := x_0 - h$, and $x_{N+1}:=x_N +h$ so that: 
$$
u''(x_0) \approx \frac{u(x_{-1}) - 2u(x_0) + u(x_{1})}{h^2} = 0
$$ 
and 
$$
u''(x_N) \approx \frac{u(x_{N-1}) - 2u(x_N) + u(x_{N+1})}{h^2} = 0.
$$
Using equations \eqref{4m1}, \eqref{4m2}, \eqref{4m4}, and \eqref{4m5} we get the following approximation
\begin{equation}
\label{4m3}
\begin{bmatrix} 
u''(x_1)\\ u''(x_2)\\ \vdots \\u''(x_{N-2}) \\u''(x_{N-1})
\end{bmatrix}
\approx 
\frac{1}{h^2}
\begin{bmatrix} 
-1 & 1 & 0&\dots & 0&0\\ 
1 & -2 &1 &\dots & 0 &0\\ 
\vdots & \ddots &\ddots & \ddots & \vdots & \vdots\\
\vdots & \vdots &\ddots & \ddots & \ddots &\vdots\\
0 & 0 &\dots &1 &-2 & 1\\ 
0 &0 &\dots & 0&-1 & 1
\end{bmatrix}
\begin{bmatrix} 
u(x_1)\\ u(x_2)\\ \vdots \\u(x_{N-2}) \\u(x_{N-1})
\end{bmatrix}.
\end{equation}

Let $v(x):= u_{xx}/(u_{xx}^2 + \epsilon)^p$. Then $v(a) = v(b) = 0$. Using similar derivations as above,
one obtains
\begin{equation}
\label{9m1}
\begin{bmatrix} 
v''(x_1)\\ v''(x_2)\\ \vdots \\v''(x_{N-2}) \\v''(x_{N-1})
\end{bmatrix}
\approx 
\frac{1}{h^2}
\begin{bmatrix} 
-2 & 1 & 0&\dots & 0&0\\ 
1 & -2 &1 &\dots & 0 &0\\ 
\vdots & \ddots &\ddots & \ddots & \vdots & \vdots\\
\vdots & \vdots &\ddots & \ddots & \ddots &\vdots\\
0 & 0 &\dots &1 &-2 & 1\\ 
0 &0 &\dots & 0&-2 & 1
\end{bmatrix}
\begin{bmatrix} 
v(x_1)\\ v(x_2)\\ \vdots \\v(x_{N-2}) \\v(x_{N-1})
\end{bmatrix}.
\end{equation}
Let
$$
\bm{D}_0:= \frac{1}{h^2}
\begin{bmatrix} 
-1 & 1 & 0&\dots & 0&0\\ 
1 & -2 &1 &\dots & 0 &0\\ 
\vdots & \ddots &\ddots & \ddots & \vdots & \vdots\\
\vdots & \vdots &\ddots & \ddots & \ddots &\vdots\\
0 & 0 &\dots &1 &-2 & 1\\ 
0 &0 &\dots & 0&-1 & 1
\end{bmatrix},\qquad 
\bm{D}_1:= 
\begin{bmatrix} 
-2 & 1 & 0&\dots & 0&0\\ 
1 & -2 &1 &\dots & 0 &0\\ 
\vdots & \ddots &\ddots & \ddots & \vdots & \vdots\\
\vdots & \vdots &\ddots & \ddots & \ddots &\vdots\\
0 & 0 &\dots &1 &-2 & 1\\ 
0 &0 &\dots & 0&-2 & 1
\end{bmatrix}
$$
and 
$$
\bm{u}: = \begin{bmatrix}u(x_1) & u(x_2) & \dots & u(x_{N-1})\end{bmatrix}^T, \qquad 
\bm{u}_0 : = \begin{bmatrix}u_0(x_1) & u_0(x_2) & \dots & u_0(x_{N-1})\end{bmatrix}^T.
$$
Then equation \eqref{25a1} leads to
$$
0 = \bm{D}_1F(\bm{u}) + \lambda(\bm{u} - \bm{u}_0),\qquad F(\bm{u}) = [\diag(\bm{D}_0\bm{u}\odot \bm{D}_0\bm{u}) + \epsilon I]^{-p}\bm{D}_0\bm{u}.
$$
Here $\odot$ denotes the Hadamard product, i.e., the element-wise matrix product, and 
the $\alpha$-th power of a diagonal matrix $\diag(d_1,...,d_n)$ where $d_i>0$ is defined as follows 
$$
\bigg(\diag(d_1,...,d_n)\bigg)^\alpha:=\diag(d_1^\alpha,...,d_n^\alpha),\qquad \alpha \in \mathbb{R}.
$$

Using similar notations as above, equations \eqref{25a3} and \eqref{25a4} become
\begin{equation}
\label{9m1}
\frac{d \bm{u}}{dt} = -\bm{D}_1F(\bm{u}) - \lambda(\bm{u} - \bm{u}_0),\qquad F(\bm{u}) = [\diag(\bm{D}_0\bm{u}\odot \bm{D}_0\bm{u}) + \epsilon I]^{-p}\bm{D}_0\bm{u}.
\end{equation}

To solve equation \eqref{9m1}, one can use the Euler method (see, e.g. \cite{Butcher}) with a step-size $h>0$ to discretize it as follows
\begin{equation}
\label{9m2}
\bm{u}_{n+1} = \bm{u}_n - h\bigg(\bm{D}_1F(\bm{u}_n) + \lambda(\bm{u}_n - \bm{u}_0)\bigg),\qquad n=1,2,... .
\end{equation}
Then $\bm{u}_n$ with sufficiently large $n$ is used as the restored data for $\bm{u}_0$. 

Again, the parameter $\lambda$ can also be chosen adaptively, based on the idea from \cite{ROF}, if $\|u_0 - f\|$ is known. Specifically, if $\|u_0 - f\| = \delta$, one would aim to use $\lambda = \lambda_n$ such that $\|\bm{u}_0 - \bm{u}_{n+1}\| = \delta$. Note that for sufficiently large $n$, i.e., at equilibrium, we would have $\bm{u}_{n+1} = \bm{u}_n$.

\subsection{2-D noise removal}
When $u_0$ is a noisy 2-D data, we propose solving the following equation and using the solution as the restored data for $u_0$:
\begin{align}
\label{x10}
0 &= \triangle \bigg(\frac{\triangle u}{(|\triangle u|^2 + \epsilon)^p}\bigg) + \lambda (u - u_0)\quad \text{in}\quad D,\\
\label{x20}
\frac{\partial u}{\partial n} &= 0\quad \text{on}\quad \partial D,\qquad \triangle u = 0\quad \text{on}\quad \partial D.
\end{align}
Here, $D\subset \mathbb{R}^2$, $\triangle$ denotes the Laplacian, and $\frac{\partial u}{\partial n}$ denotes normal derivative of $u$.

To find an approximate solution to equations \eqref{x10} and \eqref{x20}, we propose using $u(t)$ for sufficiently large $t$ as the approximate solution, where $u(t)$ solves the following evolution equation:
\begin{align}
\label{x1}
\frac{\partial u}{\partial t} &= -\triangle \bigg(\frac{\triangle u}{(|\triangle u|^2 + \epsilon)^p}\bigg) - \lambda (u - u_0)\quad \text{in}\quad D,\\
\label{x2}
\frac{\partial u}{\partial n} &= 0\quad \text{on}\quad \partial D,\qquad \triangle u = 0\quad \text{on}\quad \partial D. 
\end{align}
When a `good' approximation for $u_0$ is available, we propose using it as the initial data $u(0,x)$, $x\in D$,
for equations \eqref{x1} and \eqref{x2}. If such a `good' approximation is not available, we suggest using $u(0,x) = u_0(x)$, $x\in D$. 

In practice, equations \eqref{x1} and \eqref{x2} can be discretized by the Euler method as demonstrated in the 1-D case above (cf \eqref{9m2}). 

For discrete approximations to the Laplacian and normal derivatives $\frac{\partial u}{\partial n}$, we adopt the approach outlined in \cite{DH}. 

\section{Numerical experiments}

\subsection{1-D noise removing}

To test our method for 1-D data, we use the following functions
$$
f(x) = \sin(2\pi x),\qquad 0\le x\le 1,
$$
and 
$$
g(x) = f(x) + \sign(x-0.2) - \sign(x-0.4)+\sign(x-0.6)-\sign(x-0.8),\qquad 0\le x\le 1.
$$
Note that $f(x)$ is a smooth function, while $g(x)$ is a piecewise smooth function (see Figure \ref{figure1} below). 

The noisy versions of $f$ and $g$ are created as follows
$$
f_\delta(x) = f(x) + \delta_{rel}\frac{error}{\|error\|}\|f\|
$$
and 
$$
g_\delta(x) = g(x) + \delta_{rel}\frac{error}{\|error\|}\|g\|
$$
where $error$ is a random function and $\|f\|$ denotes the norm of the function $f$. 

In our experiments, we discretize the interval [0,1] using equidistance nodal points $(x_i)_{i=0}^n$, $x_i = i/n, i=0,...,n$. The error values at $(x_i)_{i=0}^n$ are generated randomly by the Matlab function $randn(n+1)$. 

Figure \ref{figure1} plots the original functions $f(x)$ and $g(x)$, along with their noisy counterparts $f_\delta$ and $g_\delta$ for $\delta_{rel} = 0.09$. This indicates that the noise level is approximately 9\% of the true functions, which is quite significant. The function $f(x)$ is smooth while the function $g(x)$ is piecewise smooth. We aim to evaluate how well the new filter performs on both smooth and piecewise smooth functions.

\begin{figure}[!h!t!b!]
\centerline{
\begin{tabular}{c}
\includegraphics[scale=0.65]{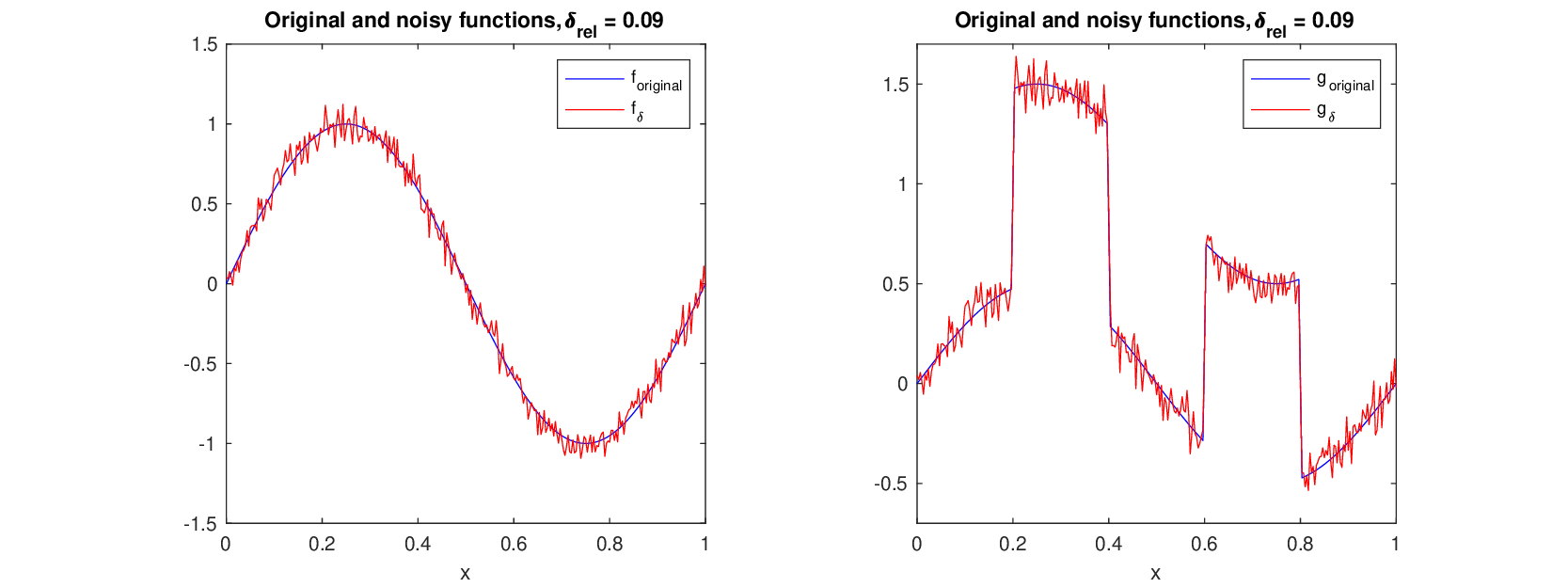}
\end{tabular}
}
\caption{\it Original and noisy functions.
}
\label{figure1}
\end{figure}

\begin{figure}[!h!t!b!]
\centerline{
\begin{tabular}{c}
\includegraphics[scale=0.65]{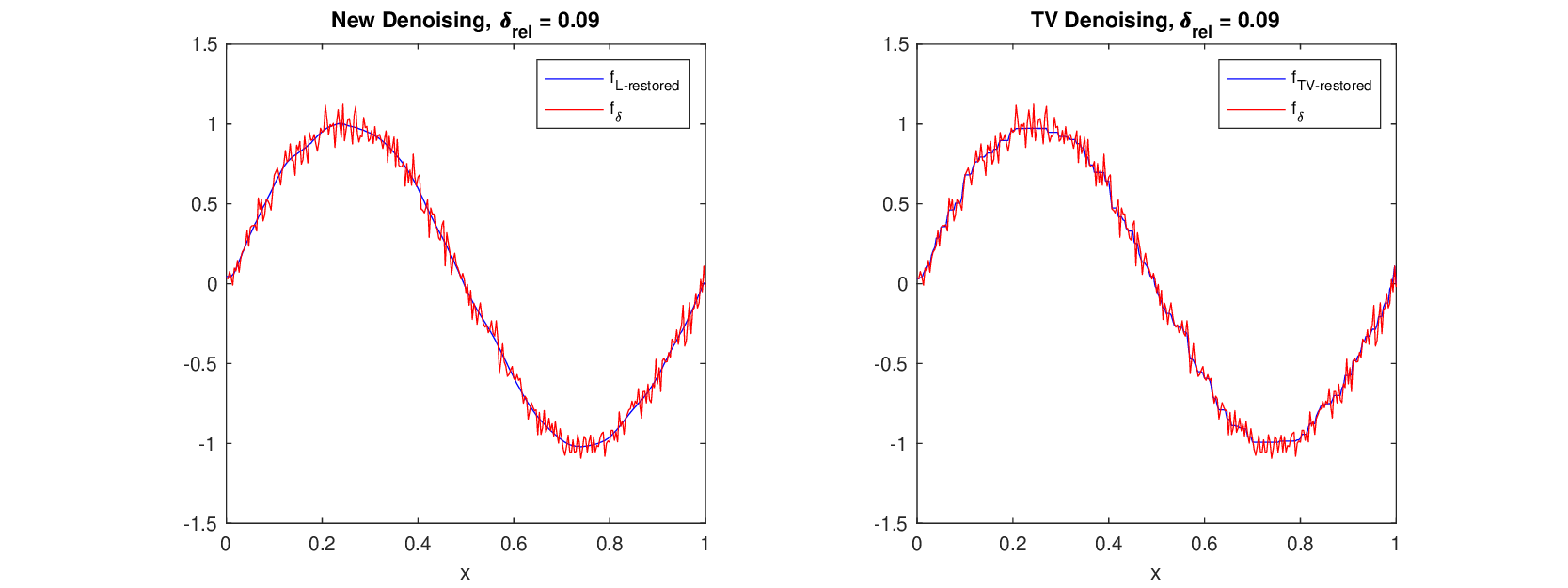}
\end{tabular}
}
\caption{\it Noisy and restored functions by the new method (left) and the TV method (right).
}
\label{figure2}
\end{figure}

Figure \ref{figure2} presents noisy function $f_\delta$ in red and the restored functions in blue using the new method (left) and the TV method (right). The figure shows that both methods effectively reduce noise from the noisy data. However, the restored function using the TV method exhibits the well-known staircase effect, appearing as a piecewise constant function. In contrast, the restored function using the new method appears to be improved compared to the TV method in this experiment.

\begin{figure}[!h!t!b!]
\centerline{
\begin{tabular}{c}
\includegraphics[scale=0.65]{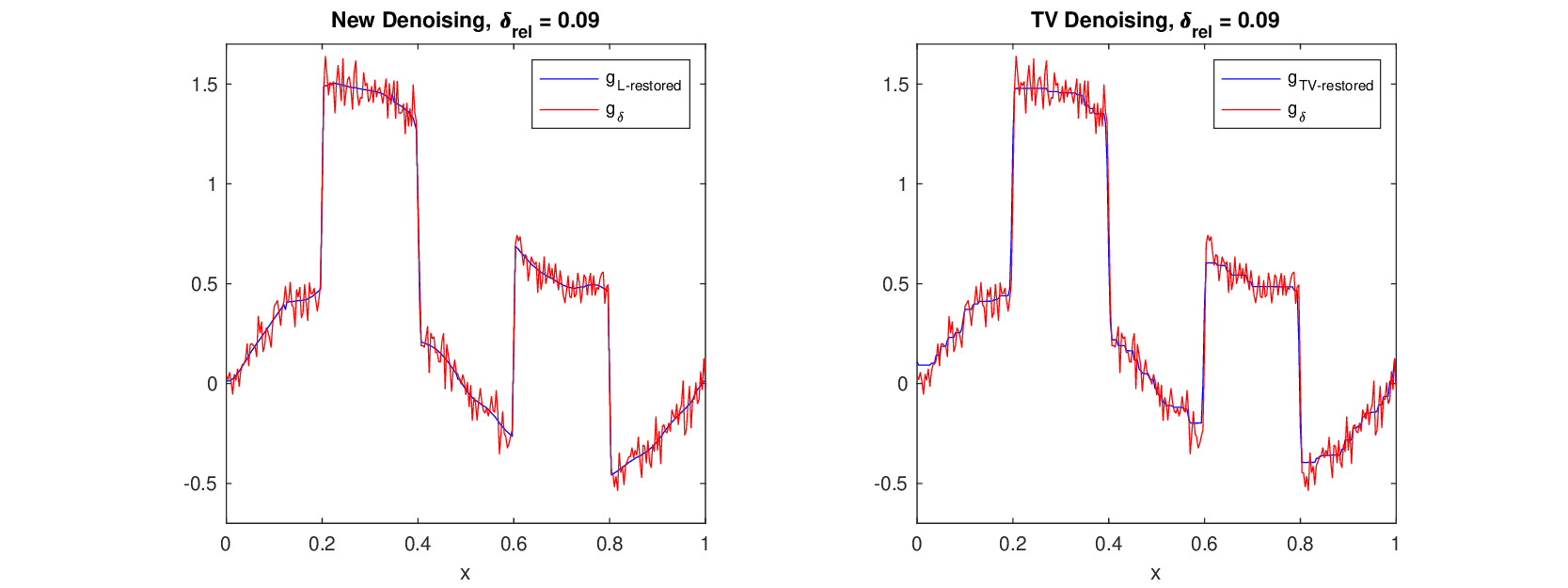}
\end{tabular}
}
\caption{\it Noisy and restored functions by the new method (left) and the TV method (right).
}
\label{figure3}
\end{figure}

Figure \ref{figure3} displays the noisy function $g_\delta$ in red and its denoised version in blue using the new method (left) and the TV method (right).  Both methods effectively reduce a significant amount of noise from $g_\delta$. 
As expected, the TV method shows the staircase effect in the reconstruction. Similar to the previous experiment, the reconstruction using the new method appears smoother and better preserves discontinuities compared to the TV method.

The conclusion from these experiments is that the new method performs very well in removing noise from 1-D piecewise smooth functions. The reconstructions using the new method are piecewise linear functions that avoid the staircase effect commonly observed with the TV method.

\subsection{2-D noise removing}

In this section, we will conduct our experiments with the following function:
$$
f(x,y) = x\sin(\pi y), \qquad (x,y)\in D = [-1,1]\times [-1,1].
$$
Using this original function, we compute the values $f_{ij}:=f(x_i,y_j)$ where
$$
x_i = -1+\frac{2i}{n-1},\qquad y_j = -1+\frac{2j}{n-1},\qquad i,j = 0,...,n-1.
$$
We assume that the matrix $U_0:=[f_{ij}]_{i,j=0}^{n-1}$ is not known, but only its noisy data $U^\delta_0:=[f^\delta_{ij}]_{i,j=0}^{n-1}$ are available. The $n\times n$ error matrix $Err_n:=[err_{ij}]_{i,j = 0}^{n-1}$ where the entries $err_{ij}:=f^\delta_{ij} - f_{ij}$ are generated by the Matlab $randn(n)$ function. The relative error quantity $\delta_{rel}$ used in our experiments is computed using the following formula: 
$$
\delta_{rel}: = \frac{\|Err_n\|}{\|U_0\|}. 
$$
Thus, we have
$$
U^\delta_0 = U_0 + \delta_{rel}\frac{Err_n}{\|Err_n\|} \|U_0\|. 
$$

Again, we use the Euler method (cf. \cite{Butcher}) to discretize equations \eqref{x1} and \eqref{x2} (see also \eqref{9m1},\eqref{9m2}), with discretizations of the Laplacian and normal derivatives adapted from \cite{DH} to solve these equations. 

\begin{figure}[!h!t!b!]
\centerline{
\begin{tabular}{c}
\includegraphics[scale=0.65]{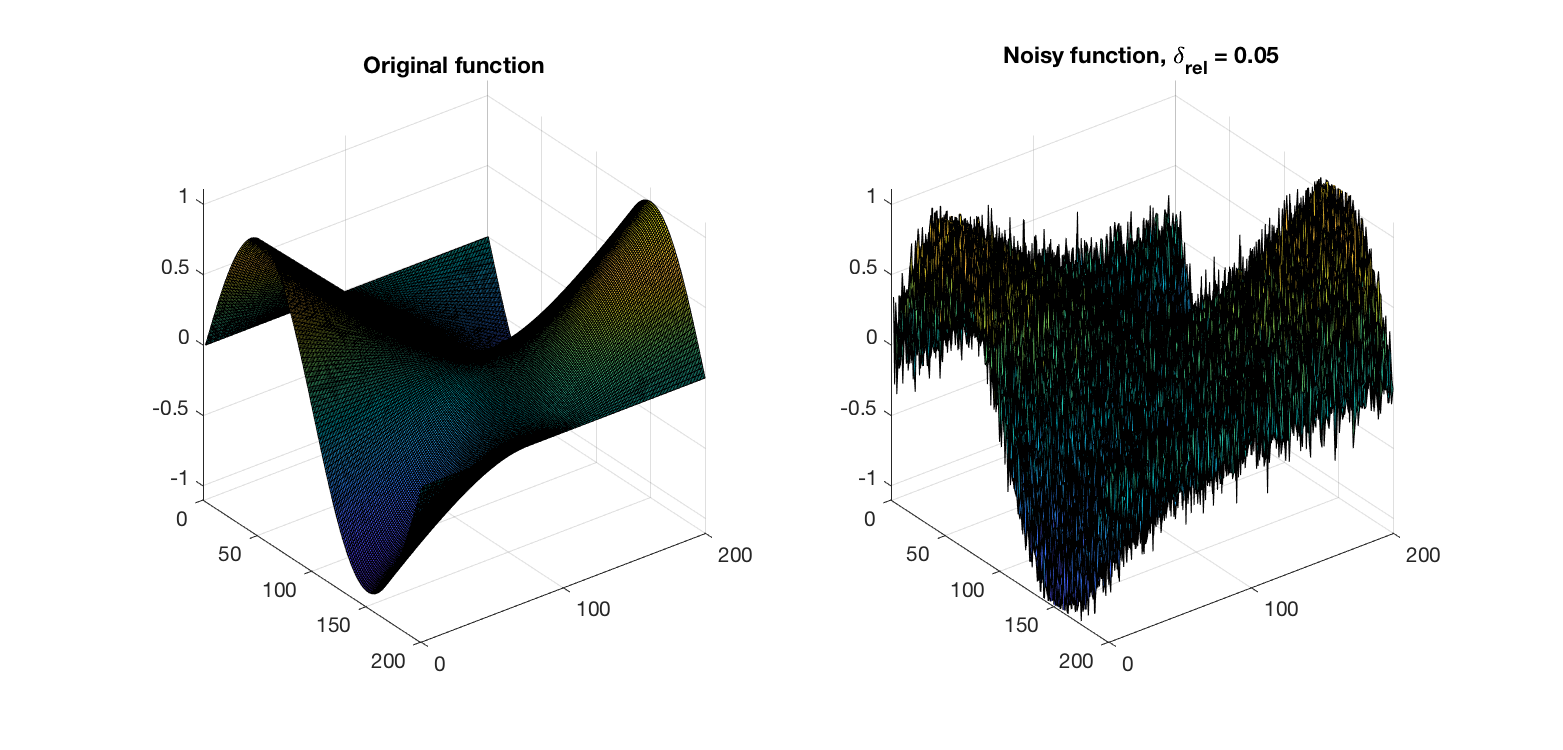}
\end{tabular}
}
\caption{\it Original and noisy functions.
}
\label{figure4}
\end{figure}

Figure \ref{figure4} plots the exact data $U_0$ and $U^\delta_0$ for $n=200$ and $\delta_{rel} = 0.05$. 
The figure shows that the noisy function exhibits a significant level of noise and looks quite different from the original function.

\begin{figure}[!h!t!b!]
\centerline{
\begin{tabular}{c}
\includegraphics[scale=0.65]{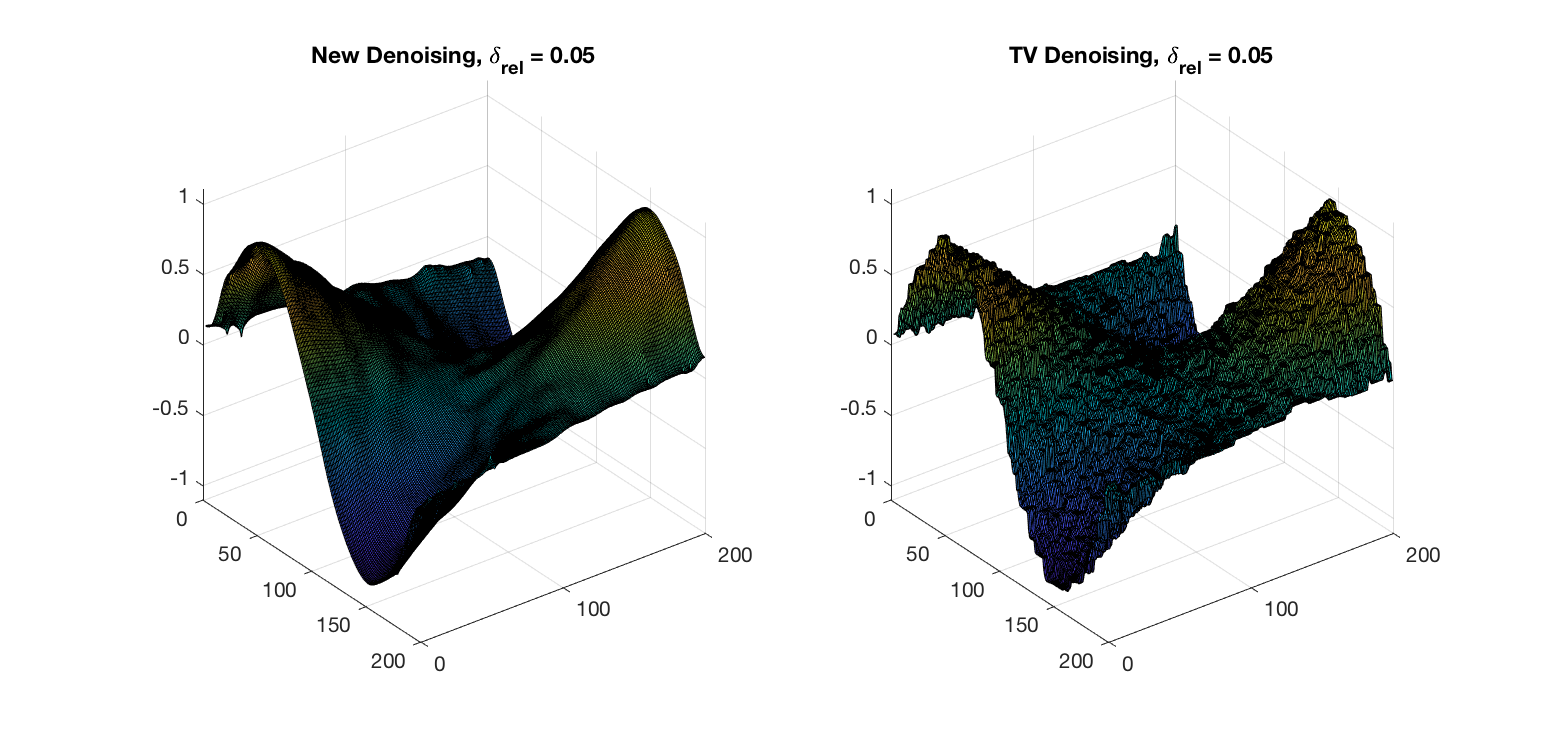}
\end{tabular}
}
\caption{\it Restored functions by the new method (left) and the TV method (right).
}
\label{figure5}
\end{figure}

Figure \ref{figure5} demonstrates the restored functions using our new method (left) and the TV method (right). As expected, our method produces a restored function that is much smoother compared to the one restored by the TV method. The TV method's restored function exhibits the usual staircase effect.

\begin{figure}[!h!t!b!]
\centerline{
\begin{tabular}{c}
\includegraphics[scale=0.65]{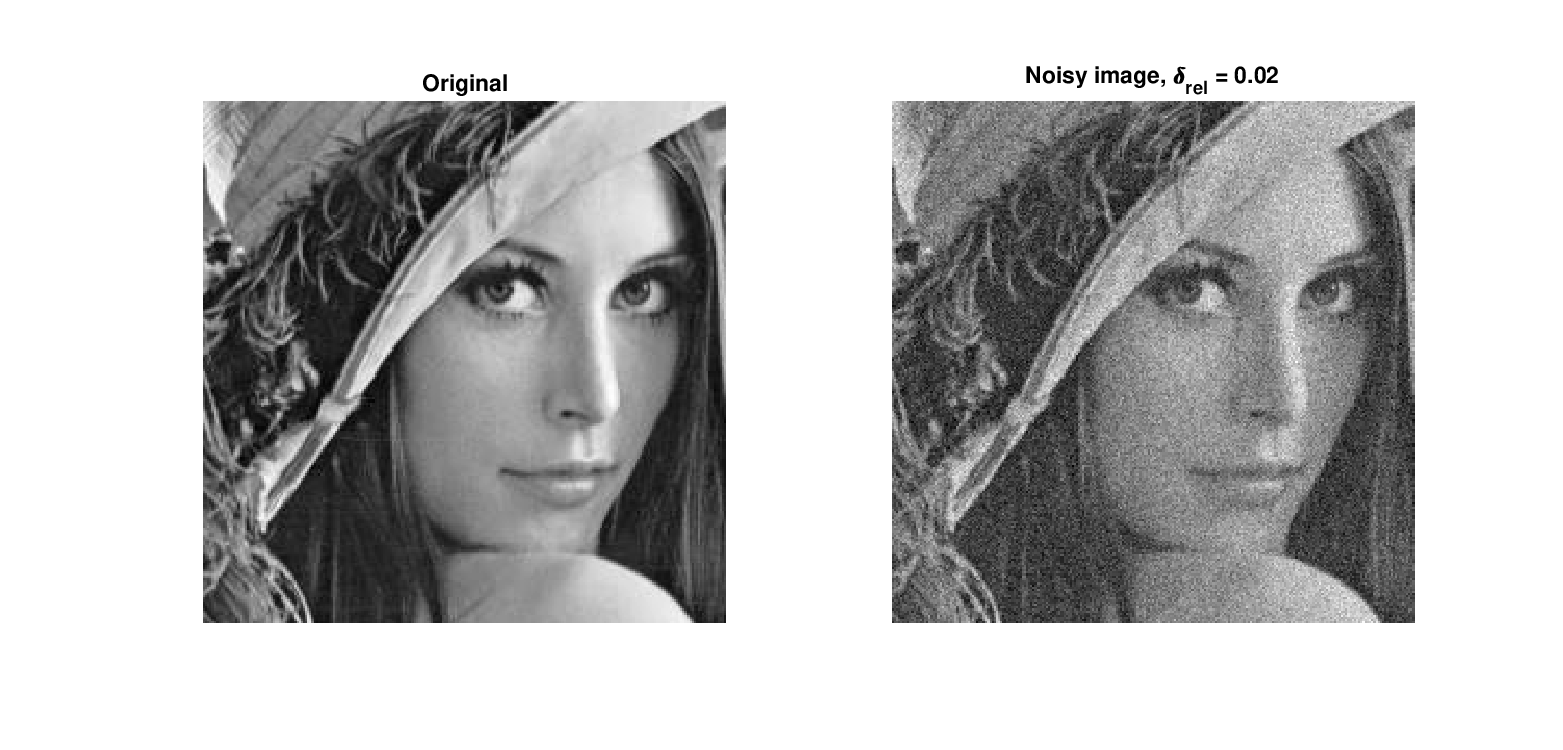}
\end{tabular}
}
\caption{\it Original and noisy images.
}
\label{figure6}
\end{figure}

We also conducted an experiment for noise removal in image restoration 
using our new method for 2D data. In particular, Figure \ref{figure6} shows an original image (left) and a noisy image (right). The noisy image is created by adding random noise to the original image. The relative error $\delta_{rel}$ 
is defined as the ratio of the 2-norm of the noise to the 2-norm of the original image.

\begin{figure}[!h!t!b!]
\centerline{
\begin{tabular}{c}
\includegraphics[scale=0.65]{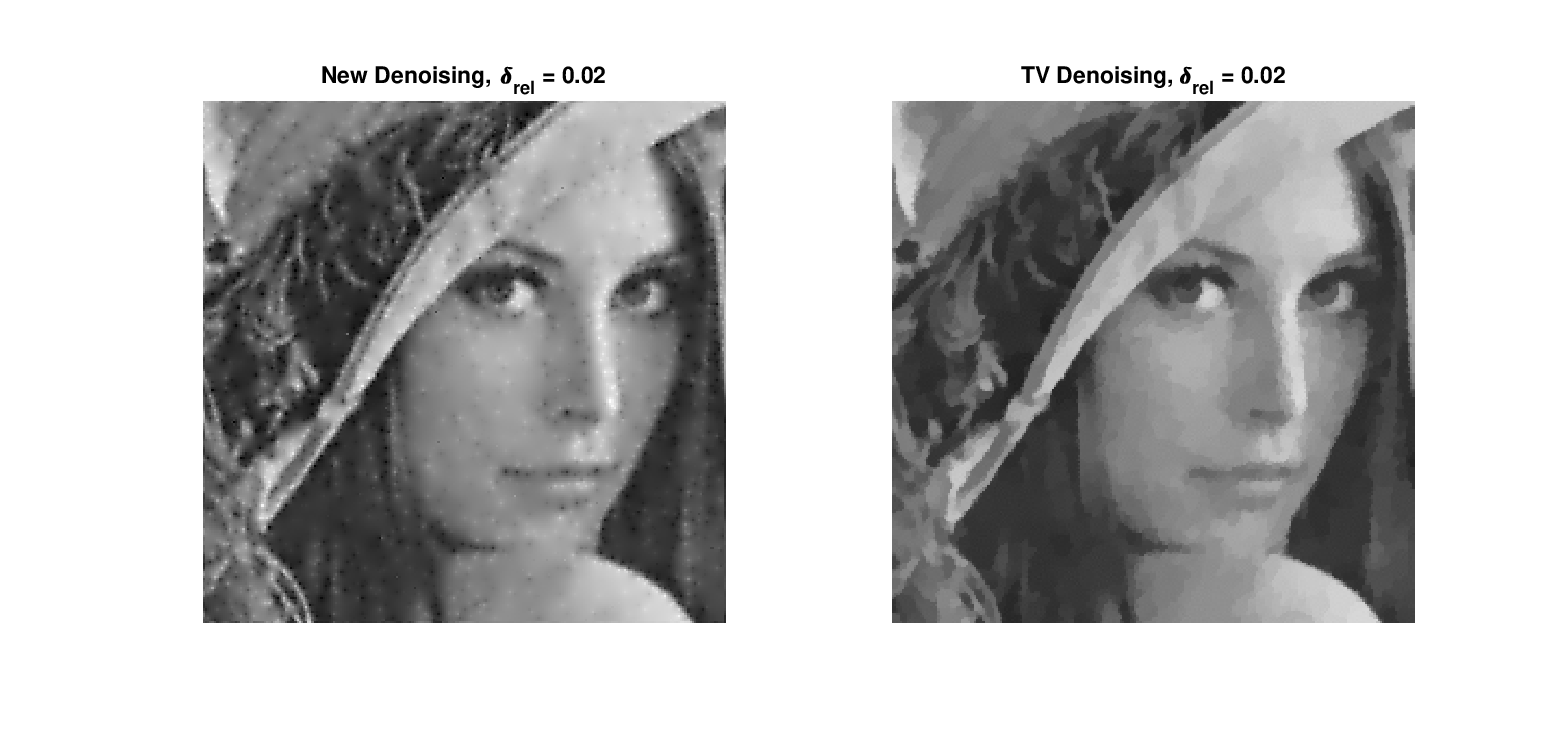}
\end{tabular}
}
\caption{\it Restored images by the new method (left) and the TV method (right).
}
\label{figure7}
\end{figure}

Figure \ref{figure7} presents the restored images using the new method (left) and the TV method (right). As observed, the restored image by the TV method exhibits the staircase effect, while the restored image by the new method does not.

\section{Conclusions}

We have proposed nonlinear filters for noise removal based on the Laplacian and discussed their numerical implementation. Our experiments have shown that the new method significantly reduces noise while preserving edges and discontinuities. Unlike the TV method, the new method does not exhibit the staircase effect commonly seen in its reconstructions.

\section{Conflicts of Interest}

The authors declare that there are no conflicts of interest regarding the publication of this paper.



\begin{thebibliography}{99}

\bibitem{BKP}
K. Bredies, K. Kunisch, T. Pock, Total generalized variation, SIAM J. Imag. Sci. 3 (3) (2010) 492--526.

\bibitem{BCM} 
P. Blomgren, T.F. Chan, P. Mulet, C. Wong, Total variation image restoration: numerical methods and extensions, in: Proceedings of the 1997 IEEE
International Conference on Image Processing, III (1997), pp. 384--387.

\bibitem{Butcher}
J. C. Butcher, Numerical Methods for Ordinary Differential Equations, Wiley, 3 edition, 2016.

\bibitem{CMM}
T. Chan, A. Marquina, P. Mulet, High-order total variation-based image restoration, SIAM J. Sci. Comput. 22 (2) (2000) 503--516.

\bibitem{DH}
S. B. Damelin and N. S. Hoang, On Surface Completion and Image Inpainting by Biharmonic Functions: Numerical Aspects, Int. J. Math. Math. Sci., vol. 2018, Article ID 3950312.

\bibitem{GQSW} 
Z. Guo, L. Qiang, J. Sun, B. Wu, Reaction-diffusion systems with growth for image denoising, Nonlinear Analysis Real World Applications 12(5) 2904--2918.


\bibitem{ROF} 
L. Rudin, S. Osher, E. Fatemi, Nonlinear total variation based noise removal algorithms, Physica D 60 (1992) 259--268.


\bibitem{SC}
L. Sun, K. Chen, A new iterative algorithm for mean curvature-based
variational image denoising, BIT Numer. Math. 54 (2), (2014) 523--553.

\bibitem{STC}
D.M. Strongand, T.F. Chan, Spatially and scale adaptive total variation based regularization and anisotropic diffusion, in: Image Processing, CAM96-
46(UCLA).



\bibitem{ZC}
J. Zhang, K. Chen, A total fractional-order variation model for image restoration
with nonhomogeneous boundary conditions and its numerical solution, SIAM
J. Imag. Sci. 8 (4) (2015) 2487--2518.


\end{thebibliography}
\end{document}